                \newcommand{\be}{\beta}
         \newcommand{\vphi}{\varphi}

\newcommand{\lb}{\lambda}

\newcommand{\cal}{\mathcal}
                      
           \newcommand{\calf}{{\cal F}}           
           \newcommand{\calp}{{\cal P}}
\newcommand{\calr}{{\cal R}}           \newcommand{\cals}{{\cal S}}
           \newcommand{\calv}{{\cal V}}

\newcommand{\arctg}{{\rm arctg}}       
       \newcommand{\Dom}{{\rm Dom}}

\newcommand{\incl}{\subseteq}          
\newcommand{\es}{\emptyset}            \newcommand{\sm}{\setminus}
\newcommand{\limpl}{\Longrightarrow}   \newcommand{\lequi}{\Longleftrightarrow}
       
\newcommand{\oo}{\infty}

             \newcommand{\sk}{\smallskip}
\newcommand{\n}{\noindent}

              \def\Rpoo{R\cup \{\oo\}}             
\def\Rmpoo{R\cup\{-\oo,\oo\}}          

\def\dtends  {\stackrel {\it d}{\longrightarrow}}

\def\Dtends  {\stackrel {\it D}{\longrightarrow}}
\def\(V)tends  {\stackrel {(\calv)}{\longrightarrow}}

\newcommand{\barr}{\begin{array}}         \newcommand{\earr}{\end{array}}
\newcommand{\bcor}{\begin{corollary}}     \newcommand{\ecor}{\end{corollary}}
\newcommand{\beq}{\begin{equation}}       \newcommand{\eeq}{\end{equation}}
\newcommand{\bit}{\begin{itemize}}        \newcommand{\eit}{\end{itemize}}
\newcommand{\blemma}{\begin{lemma}}       \newcommand{\elemma}{\end{lemma}}
\newcommand{\bproof}{\begin{proof}}       \newcommand{\eproof}{\end{proof}}
\newcommand{\bprop}{\begin{proposition}}  \newcommand{\eprop}{\end{proposition}}
\newcommand{\brem}{\begin{remark}}        \newcommand{\erem}{\end{remark}}
\newcommand{\btab}{\begin{tabular}}       \newcommand{\etab}{\end{tabular}}
\newcommand{\btheorem}{\begin{theorem}}   \newcommand{\etheorem}{\end{theorem}}

\documentclass[reqno]{amsart}

\usepackage{amssymb}
\usepackage{latexsym}
\usepackage{amsfonts,euscript}

\newtheorem{theorem}{\bf Theorem}
\newtheorem{corollary}{\bf Corollary}

\newtheorem{lemma}{\bf Lemma}
\newtheorem{proposition}{\bf Proposition}
\newtheorem{remark}{\bf Remark}

\begin{document}

\title
[Gauge Brezis-Browder Principles and Dependent Choice]
{GAUGE BREZIS-BROWDER PRINCIPLES \\
AND DEPENDENT CHOICE}

\author{Mihai Turinici}
\address{
"A. Myller" Mathematical Seminar;
"A. I. Cuza" University;
700506 Ia\c{s}i, Romania
}
\email{mturi@uaic.ro}


\subjclass[2010]{
49J53 (Primary), 47J30 (Secondary).
}

\keywords{
Quasi-order, maximal element, gauge structure, Dependent Choice Principle, 
lsc function, inf-lattice, Lipschitz map, discrete space. 
}

\begin{abstract}
The gauge Brezis-Browder Principle in
Turinici [Bull. Acad. Pol. Sci. (Math.), 30 (1982), 161-166]
is obtainable from the 
Principle of Dependent Choices (DC)
and implies  
Ekeland's Variational Principle
(EVP);
hence, it is equivalent with both (DC) and (EVP).
This is also true for 
the gauge variational principle deductible from it, 
including the one in Bae, Cho, and Kim
[Bull. Korean Math. Soc. 48 (2011), 1023-1032].
\end{abstract}

\maketitle

\section{Introduction}
\setcounter{equation}{0}

Let $M$ be a nonempty set. Take a {\it quasi-order} $(\le)$
(i.e.: reflexive and transitive relation) over it;
and a function $\vphi: M\to \Rmpoo$.
Call the point $z\in M$, $(\le,\vphi)$-{\it maximal} when:
$z\le w\in M$  implies $\vphi(z)=\vphi(w)$;
or, equivalently: 
$\vphi$ is constant on $M(z,\le):=\{x\in M; z\le x\}$;
the set of all these will be denoted as $\max(M;\le;\vphi)$.
A basic result about  such points is the 1976
Brezis-Browder ordering principle
\cite{brezis-browder-1976} (in short: BB).

\btheorem \label{t1}
Suppose that
\bit
\item[(a01)]  
$(M,\le)$ is sequentially inductive: \\
each ascending sequence has an upper bound (modulo $(\le)$)
\item[(a02)]
$\vphi$ is $(\le)$-decreasing\ 
($x_1\le x_2$ $\limpl$ $\vphi(x_1)\ge \vphi(x_2)$) 
\item[(a03)]
$\vphi(M)\incl R$ and\ $\vphi$  
is bounded below ($\inf \vphi(M)> -\oo$).
\eit
Then, $\max(M;\le;\vphi)$ is 

{\bf i)} $(\le)$-cofinal in $M$
[for each $u\in M$ there exists $v\in \max(M;\le;\vphi)$ with $u\le v$]

{\bf ii)} $(\le)$-invariant in $M$
[$z\in \max(M;\le;\vphi)$ $\limpl$ $M(z,\le) \incl \max(M;\le;\vphi)$].
\etheorem

This statement includes (cf. Section 4)
Ekeland's Variational Principle \cite{ekeland-1979} (in short: EVP); 
and found some useful applications to convex and non-convex analysis
(see the above references). So, it was the subject of
many extensions; see, for instance,
Hyers, Isac and Rassias \cite[Ch 5]{hyers-isac-rassias-1997}.
These are interesting from a technical perspective;
but, in all concrete situations when a
variational principle of this type (VP, say) is to be applied, 
a substitution of it by the Brezis-Browder's is always possible.
On the other hand (cf. Section 3), any VP like before  
is reducible to the 
Bernays-Tarski  Dependent Choice Principle
(in short: DC), discussed in Section 2.
This ultimately raises the question of to what extent 
are the inclusions  
(DC) $\limpl$ (BB) $\limpl$ (EVP) effective.
A negative answer to this is to be found in 
Turinici \cite{turinici-2011}.
Here, we shall be concerned with the
inclusion between their gauge versions.
Precisely, we show in Section 4  that
(DC)  $\limpl$ (BBg) $\limpl$ (EVPg) $\limpl$ (EVP); 
here, (BBg) is the gauge  variant of (BB) in
Turinici \cite{turinici-1982}
and (EVPg) is the gauge version of (EVP).
This, along with (EVP) $\limpl$ (DC) (cf. Section 5),
closes the circle between all these.
In particular, the gauge 
variational principle in
Bae, Cho, and Kim \cite{bae-cho-kim-2011}
enters in such a chain.
Further aspects will be delineated elsewhere.

\section{Dependent Choice Principles}
\setcounter{equation}{0}

Let $M$ be a nonempty set.
By a {\it relation} over it we mean any  
mapping $\calr$  from $M$ to $\calp(M)$ 
(=the class of all subsets in $M$). 
As usually, we identify $\calr$ with its graph in $M\times M$;
so, given $x,y\in M$, we may write $y\in \calr(x)$ as $x\calr y$. 
Call the relation $\calr$, {\it proper} when 
\bit
\item[(b01)]
$\calr(x)$ is nonempty, for each $x\in M$.
\eit
Note that, under such a condition, $\calr$ acts as a mapping 
between $M$ and $\calp_0(M)$ 
(=the subclass of all nonempty parts in $M$).

{\bf (A)}
The following "Dependent Choice" principle (in short: DC) 
is our starting point.
Given $a\in M$, call the sequence $(x_n; n\ge 0)$ in $M$,
{\it $(a;\calr)$-iterative} provided
\bit
\item[(b02)]
$x_0=a$,\ $x_{n+1}\in \calr(x_n)$,\  for all $n\ge 0$.
\eit

\bprop \label{p1}
Let $\calr$ be a proper relation over $M$.
Then, for each $a\in M$ there exists 
at least one $(a,\calr)$-iterative sequence in $M$.
\eprop

This statement -- due, independently, to
Bernays \cite{bernays-1942}
and
Tarski \cite{tarski-1948} --
has a strong connection with 
the Axiom of Choice (in short: AC) 
from the usual Zermelo-Fraenkel axiomatic system (ZF),
as described in Cohen \cite[Ch 2, Sect 3]{cohen-1966}.
Precisely, in the {\it reduced} 
Zermelo-Fraenkel system (ZF-AC),
we have (AC) $\limpl$ (DC); but not conversely;
see, for instance, 
Wolk \cite{wolk-1983}.
Moreover, the {\it (DC)-added Zermelo-Fraenkel} system
(ZF-AC+DC) is large enough so as to 
include the "usual" mathematics; 
see, for instance, 
Moskhovakis \cite[Ch 8]{moskhovakis-2006}.

{\bf (B)}
Let $(\calr_n; n\ge 0)$ be a sequence of relations on $M$.
Given $a\in M$, let us say that the sequence $(x_n; n\ge 0)$ in $M$
is {\it $(a;(\calr_n; n\ge 0))$-iterative} provided
\bit
\item[(b03)]
$x_0=a$, $x_{n+1}\in \calr_n(x_n)$, $\forall n$.
\eit
The following 
"Diagonal Dependent Choice" principle
(in short: DDC) is available.

\bprop \label{p2}
Let $(\calr_n; n\ge 0)$ be a sequence of proper relations on $M$.
Then, for each $a\in M$ there exists 
at least one $(a; (\calr_n; n\ge 0))$-iterative sequence in $M$.
\eprop

Clearly, (DDC) includes (DC); to which it reduces when
$(\calr_n; n\ge 0)$ is constant.
The reciprocal of this is also true.
In fact, letting the premises of (DDC) hold,
put $P=N\times M$; and let $\cals$ be the relation over $P$ 
introduced as
\bit
\item[]
$\cals(i,x)=\{i+1\}\times \calr_i(x)$,\ \ $(i,x)\in P$.
\eit
It will suffice applying (DC) to $(P,\cals)$ and $b:=(0,a)\in P$
to get the conclusion in the statement; we do not give details.

Summing up, (DDC) is provable in (ZF-AC+DC).
This is valid as well for its variant,
referred to as: the "Selected Dependent Choice" principle 
(in short: SDC).

\bprop \label{p3}
Let the map $F:N\to \calp_0(M)$ and the 
relation $\calr$ over $M$ fulfill
\bit
\item[(b04)]
($\forall n\in N$):\ $\calr(x)\cap F(n+1)\ne \es$,\ \ 
$\forall x\in F(n)$.
\eit
Then, for each $a\in F(0)$ there exists 
a sequence $(x(n); n\ge 0)$ in $M$ with  
\beq \label{201}
x(0)=a;\ x(n)\in F(n),\ \forall n;\ x(n+1)\in \calr(x(n)),\ \forall n.
\eeq
\eprop

As before, (SDC) $\limpl$ (DC) ($\lequi$ (DDC));
just take $F(n)=M$, $n\ge 0$.
But, the reciprocal is also true, in the sense:
(DDC) $\limpl$ (SDC). This follows from 

\bproof {\bf (Proposition \ref{p3})}
Let the premises of (SDC) be true.
Define a sequence of relations $(\calr_n; n\ge 0)$ 
over $M$ as: for each $n\ge 0$,
\bit
\item[(b05)]
$\calr_n(x)=\calr(x)\cap F(n+1)$,\ \ if $x\in F(n)$, \\
$\calr_n(x)=\{x\}$, otherwise ($x\in M\sm F(n)$).
\eit
Clearly, $\calr_n$ is proper, for all $n\ge 0$. 
So, by (DDC), it follows that, 
for the starting $a\in F(0)$, there exists 
a sequence $(x(n); n\ge 0)$ in $M$ with
the property (b03).
Combining with the very definition (b05), it 
follows that  (\ref{201}) is holding. 
\eproof

{\bf (C)}
In particular, when $\calr=M\times M$, 
(b04) holds. The corresponding variant of (SDC) is just 
(AC(N)) (=the Denumerable Axiom of Choice).
Precisely, we have

\bprop \label{p4}
Let $F:N \to \calp_0(M)$ be a function.
Then, for each $a\in F(0)$ there exists a
function $f:N\to M$ with 
$f(0)=a$ and $f(n)\in F(n)$, $\forall n\ge 0$.
\eprop

\brem \label{r1}
\rm
Note that, as a consequence of the above facts, 
(DC) $\limpl$ (AC(N)), in (ZF-AC).
A direct verification of this is obtainable by 
taking $P=N\times M$ 
and introducing the relation $\calr$ over it, according to:
\bit
\item[]
$\calr(n,x)=\{n+1\}\times F(n+1)$,\ \ $n\ge 0$, $x\in M$;
\eit
we do not give details.
The reciprocal of the written inclusion is not true;
see
Moskhovakis \cite[Ch 8, Sect 8.25]{moskhovakis-2006}
for details.
\erem

\section{Gauge ordering principles}
\setcounter{equation}{0}

Let $M$ be a nonempty set; remember that 
$\calp_0(M)=\{Y\in \calp(M); Y\ne \es\}$. 
As already specified, the axiomatic 
system in use is (ZF).
\sk

{\bf (A)}
Given some property $\pi$ involving $\calp_0(M)$,
denote by $(\pi)$ the subclass of all $Y\in \calp_0(M)$ 
fulfilling it.
In this case, let us say that $\pi$ is {\it inductive} provided:
\bit
\item[(c01)]
($Y_i\in  (\pi)$, $\forall i\ge 0$) implies $Y:=\cap\{Y_i; i\ge 0\} \in (\pi)$\ 
(hence, $Y\in \calp_0(M)$).
\eit
An interesting example of this type is the following. 
Let $(M,\le)$ be a quasi-ordered structure. 
Call $Z\in \calp_0(M)$, {\it $(\le)$-cofinal} in $M$
when [$M(u,\le)\cap Z\ne \es$, $\forall u\in M$].
In addition, let us say that 
$Z\in \calp_0(M)$ is {\it $(\le)$-invariant} 
provided $w\in Z$ implies $M(w,\le)\incl Z$.
The intersection of these properties 
will be referred to as: 
$Z$ is $(\le)$-cofinal-invariant; 
in short: $(\le)$-cof-inv.

\bprop \label{p5}
Assume that $(M,\le)$ is sequentially inductive (cf. (a01)).
Then, the $(\le)$-cof-inv property 
is inductive (in (ZF-AC+DC)).
\eprop

\bproof
Let $(F(i); i\ge 0)$
be a sequence in $\calp_0(M)$ such that: 
$F(i)$ is $(\le)$-cof-inv, for each $i\ge 0$.
We intend to show that 
$Y:=\cap\{F(i); i\ge 0\}$ is endowed 
with the same property.
Clearly, $Y$ is $(\le)$-invariant; 
but, for the moment, $Y=\es$ cannot be avoided.
We show that $Y$ is $(\le)$-cofinal too; hence nonempty.
Let $u\in M$ be arbitrary fixed.
Further, let the relation $\calr$ over $M$ be introduced as
[$\calr(x)=M(x,\le)$, $x\in M$];
i.e.: $\calr$ is the {\it graph} of $(\le)$.
By the $(\le)$-cofinal property,
\beq \label{301}
\calr(x)\cap F(i)=M(x\le)\cap F(i)\ne \es,\ 
\forall i\ge 0, \forall x\in M.
\eeq
In particular, this tells us that 
$M(u,\le)\cap F(0)\ne \es$;
let $a$ be one of its elements.
From Proposition \ref{p3} it follows that, for 
this starting element, there exists a sequence 
$(x_n; n\ge 0)$ in $M$ with
\beq \label{302}
x_0=a;\ x_n\in F(n), \forall n;\  x_n\le x_{n+1}, \forall n.
\eeq
As $(M,\le)$ is sequentially inductive, there exists 
at least one $v\in M$ with $x_n\le v$, $\forall n$.
In particular, from $u\le a=x_0\le v$, one has $u\le v$.
Moreover, by the $(\le)$-invariance properties
of our sequence, we have $v\in F(n)$, $\forall n$; 
hence $v\in Y$. The proof is complete.
\eproof

{\bf (B)}
Let again $(M,\le)$ be a quasi-ordered structure;
and $\vphi:M\to R\cup \{-\oo,\oo\}$ be a function.
Define the $(\le,\vphi)$-maximal property of some
$z\in M$  as in Section 1; 
remember that the class of all these
was denoted as $\max(M;\le;\vphi)$.
Technically speaking, 
sufficient conditions for existence of such elements 
are to be written in terms of the underlying function $\vphi$
belonging to certain subclasses (=subsets) 
of $\calf(M,R\cup \{-\oo,\oo\})$.
[Here, for each couple $A,B$ of nonempty sets, 
$\calf(A,B)$ stands for the class of all functions from 
$A$ to $B$; when $A=B$, we write $\calf(A)$ in place of 
$\calf(A,A)$].
The basic ones are listed below:
\bit
\item[(P1)] $\vphi(M)\incl R$ and $\vphi$ is bounded 
($-\oo< \inf \vphi(M)\le \sup \vphi(M)< \oo$)
\item[(P2)] 
general case\
($\vphi(M)\cap \{-\oo,\oo\} \ne \es$ cannot be avoided)
\item[(P3)] 
$\vphi(M)\incl \Rpoo$ and $\vphi$ is bounded below ($\inf \vphi(M)> -\oo$)
\item[(P4)]  
$\vphi(M)\incl \Rpoo$ and $\vphi$ is positive ($\inf \vphi(M)\ge 0$).
\item[(P5)] 
$\vphi(M)\incl R$ and\ $\vphi$  
is bounded below ($\inf \vphi(M)> -\oo$)\ \ [cf. (a03)]
\item[(P6)]  
$\vphi(M)\incl M$ and $\vphi$ 
is positive ($\inf \vphi(M)\ge 0$). 
\eit

The following "multiple" 
ordering principle is now considered:

\btheorem \label{t2}
Assume that (a01) and (a02) are valid; as well as
\bit
\item[(c02)]
$\vphi$ belongs to the subclass (Pj)\ 
(for some $j\in \{1,2,3,4,5,6\}$).
\eit
Then, $\max(M;\le;\vphi)$ is 
$(\le)$-invariant  
and $(\le)$-cofinal
(hence, nonempty) in $M$.
\etheorem

For simplicity, we shall indicate this 
ordering principle as (BB-Pj) [where $j\in \{1,2,3,4,5,6\}$].
Note that (BB-P2) is the "extended" variant of (BB) due to 
C\^{a}rj\u{a} and Ursescu \cite{carja-ursescu-1993};
referred to as the C\^{a}rj\u{a}-Ursescu variational principle;
(in short: CU).
Moreover, (BB-P5) is just (BB)
[stated in Section 1]. 
\sk

The relationships between these principles 
are clarified in

\blemma \label{le1}
We have (in (ZF-AC)):
\beq \label{303}
\mbox{
(BB-P1) $\limpl$ (BB-P2) $\limpl$ (BB-P3) $\limpl$ (BB-P5)
$\limpl$ (BB-P1)
}
\eeq
\beq \label{304}
\mbox{
(BB-P3) $\lequi$ (BB-P4),\  (BB-P5) $\lequi$ (BB-P6).
}
\eeq
\elemma

\bproof
{\bf i)}
The inclusion (BB-P4) $\limpl$ (BB-P3) and (BB-P6) $\limpl$ (BB-P5)
are deductible from the following remark: 
if the function $\vphi$ is like in (BB-P3) (resp., (BB-P5)), then
[$\psi(.)=\vphi(.)-\inf\vphi(M)$]  fulfills the requirements of
(BB-P4) (resp., (BB-P6)). 
This, along with the reciprocal inclusions being fulfilled,
proves (\ref{304}).

{\bf ii)}
The inclusions in (\ref{303}), 
with the exception of the first one are immediate.

{\bf iii)}
It remains to verify the quoted relation.
Let the premises of (BB-P2) hold.
Define the  function  $\chi:M \to [0, \pi]$ as
[$\chi(x)=A(\vphi(x)), x\in M$]; where
\bit
\item[]
$A(t)=\pi/2+\arctg(t)$ if $t\in R$; $A(-\oo)=0$; $A(\oo)=\pi$.
\eit
Clearly, $\chi$ fulfills (a02) and 
belongs to the subclass (P1). 
Therefore, by the conclusion of (BB-P1), 
for each $u\in M$ there exists a $(\le,\chi)$-maximal
$v\in M$ with $u\le v$. This, along with 
$\max(M;\le;\vphi)=\max(M;\le;\chi)$, 
gives the desired conclusion. 
\eproof

Note that, the obtained relations cannot assure us that 
these principles are deductible in (ZF-AC+DC).
This, however, holds; as results from

\bprop \label{p6}
We have (in (ZF-AC)) (DC) $\limpl$ (BB-P1);
hence (by the above) (DC) $\limpl$ (BB-Pj), 
for each $j\in \{1,2,3,4,5,6\}$.
\eprop

A complete proof of this may be found in 
Turinici \cite{turinici-2011};
see also
C\^{a}rj\u{a}, Necula and Vrabie 
\cite[Ch 2, Sect 2.1]{carja-necula-vrabie-2007}.
For completeness, we shall sketch the
argument (in our new setting).

\bproof {\bf (Proposition \ref{p6})}
The $(\le)$-invariance property of $\max(M;\le;\vphi)$
is clear; so, it remains to establish the $(\le)$-cofinal property
of the same. 
So, assume that (a01)+(a02) hold; and that $\vphi$ 
is in the subclass (P1).
Define the function $\be:M\to R$ as:
$\be(v):=\inf[\vphi(M(v,\le))]$, $v\in M$.
Clearly, $\be$ is increasing and 
\beq \label{305}
\mbox{
$\vphi(v)\ge \be(v)$, for all $v\in M$.
}
\eeq
Moreover, (a02) gives at once a characterization like
\beq \label{306}
\mbox{
$v$ is $(\le,\vphi)$-maximal iff $\vphi(v)=\be(v)$.
}
\eeq
Assume by contradiction that the $(\le)$-cofinal property
is false; i.e. [in combination with (\ref{306})]
there must be some $u\in M$ such that:
\bit
\item[(c03)]
for each $v\in M_u:=M(u,\le)$, one has $\vphi(v)> \be(v)$.
\eit
Consequently (for all such $v$), 
$\vphi(v)> (1/2)(\vphi(v)+\be(v)) > \be(v)$;
hence 
\beq \label{307}
\mbox{
$v\le w$ and $(1/2)(\vphi(v)+\be(v))> \vphi(w)$, 
}
\eeq
for at least one $w$ (belonging to $M_u$).
The relation $\calr$ over $M_u$ introduced via (\ref{307})
is then proper (cf. (b01)).
So, by (DC), there must be 
a sequence $(u_n)$ in $M_u$ with $u_0=u$ and
\beq \label{308}
\mbox{
$u_n\le u_{n+1}$, $(1/2)(\vphi(u_n)+\be(u_n))> \vphi(u_{n+1})$, 
for all $n$.
}
\eeq
We have thus constructed an ascending sequence  $(u_n)$ in $M_u$
for which the real sequence $(\vphi(u_n))$ is (by (c03)) 
strictly descending and bounded below; hence 
$\lb:=\lim_n \vphi(u_n)$ exists in $R$.
By (a01), $(u_n)$ is bounded from above in $M$; 
i.e., there exists  $v\in M$ such that $u_n\le v$, for all $n$.
From (a02), $\vphi(u_n)\ge \vphi(v)$, $\forall n$; 
and (by the properties of $\be$)
$\vphi(v)\ge \be(v)\ge \be(u_n)$, $\forall n$.
The former of these relations gives $\lb\ge \vphi(v)$. 
On the other hand, the latter of these relations yields
(via (\ref{308}))
$(1/2)(\vphi(u_n)+\be(v))> \vphi(u_{n+1})$, for all $n\in N$. 
Passing to limit as $n\to \oo$  gives
$(\vphi(v)\ge )\be(v)\ge \lb$; so, combining with the preceding one,
$\vphi(v)=\be(v)(=\lb)$, contradiction.
Hence, (c03) cannot be accepted; and the conclusion follows.
\eproof

In particular, the equivalent (in (ZF-AC))
ordering principles (BB) and (CU) are deductible 
(again in (ZF-AC)) from (DC).
For the reciprocal inclusions, we refer to Section 5 below.
\sk

{\bf (C)}
A denumerable version of these facts may be given as follows.
Let $\Phi=(\vphi_i; i\ge 0)$ be a sequence of maps in 
$\calf(M,R\cup\{-\oo,\oo\})$;
it will be referred to as a {\it gauge function} over
$\calf(M,R\cup\{-\oo,\oo\})$.
Call $z\in M$, $(\le,\Phi)$-{\it maximal}, provided  
$z$ is $(\le,\vphi_i)$-maximal, for each $i\ge 0$.
The class of all these will be denoted as 
$\max(M;\le;\Phi)$; hence, by definition, 
$\max(M;\le;\Phi)=\cap\{\max(M;\le;\vphi_i); i\ge 0\}$.
To get an existence result for such points,
let us accept, in addition to (a01),
\bit
\item[(c04)]
$\Phi$ is decreasing:\ $\vphi_i$ is decreasing, $\forall i\ge 0$.
\eit
Further, for each $j\in  \{1,2,3,4,5,6\}$, 
let [Pj] stand for the [attached to (Pj)] subclass (=subset) of 
all gauge functions class over $\calf(M,R\cup\{-\oo,\oo\})$
introduced as:
\bit
\item[(c05)]
$\Phi$ belongs to the subclass [Pj] iff 
$\vphi_i$ belongs to the subclass (Pj),\ $\forall i\ge 0$.
\eit

The following "multiple" gauge
ordering principle enters into discussion:

\btheorem \label{t3}
Assume that (a01)+(c04) are valid; as well as
\bit
\item[(c06)]
$\Phi$ belongs to the subclass [Pj]\ (for some $j\in \{1,2,3,4,5,6\}$).
\eit
Then, $\max(M;\le;\Phi)$ is $(\le)$-cofinal in $M$
[for each $u\in M$ there exists a $(\le, \Phi)$-maximal
$v\in M$ with $u\le v$]
and $(\le)$-invariant in $M$
[$u$ is $(\le,\Phi)$-maximal and $u\le v$ imply 
$v$ is $(\le,\Phi)$-maximal].
\etheorem

For simplicity, we shall indicate these gauge
ordering principles as (BBg-Pj) [where $j\in \{1,2,3,4,5,6\}$].
Note that (BBg-P2) is the gauge variant of the ordering 
principle (CU) (see above); so that, it will be 
indicated as (CUg).
On the other hand, (BBg-P5) is nothing else than the 
gauge variant of (BB) obtained in 
Turinici \cite{turinici-1982};
denoted as (BBg).
\sk

The relationships between these are clarified in

\blemma \label{le2}
We have (in (ZF-AC)):
\beq \label{309}
\mbox{
(BBg-P1) $\limpl$ (BBg-P2) $\limpl$ (BBg-P3) $\limpl$ (BBg-P5)
$\limpl$ (BBg-P1)
}
\eeq
\beq \label{310}
\mbox{
(BBg-P3) $\lequi$ (BBg-P4),\  (BBg-P5) $\lequi$ (BBg-P6).
}
\eeq
\elemma

The proof mimics that of Lemma \ref{le1}; 
so, it will be omitted.

As before, the obtained relations cannot assure us that 
the principles in question are deductible
in (ZF-AC+DC). This, however, holds; as follows from

\bprop \label{p7}
We have (in (ZF-AC)) 
(DC) $\limpl$ (BBg-Pj),\ $j\in \{1,2,3,4,5,6\}$.
\eprop

\bproof
Let the premises of (BBg-Pj) be accepted.
From (BB-Pj), we have that
$Y_i:=\max(M;\le;\vphi_i)$ is nonempty $(\le)$-cof-inv,
for each $i\ge 0$.
This, along with Proposition \ref{p5}
(valid in (ZF-AC+DC)), tells us that 
$\cap\{Y_i; i\ge 0\}=\max(M;\le;\Phi)$
has the same properties; and conclusion follows.
\eproof

\brem \label{r2}
\rm

By the very arguments above, one gets, 
in (ZF-AC+DC):
\beq \label{311}
\mbox{
(BB-Pj) $\limpl$ (BBg-Pj) $\limpl$ (BB-Pj),\ $j\in \{1,2,3,4,5,6\}$.
}
\eeq
Hence, for each $j\in \{1,2,3,4,5,6\}$, 
the ordering principle (BB-Pj) is equivalent 
with its gauge version (BBg-Pj).
This, however, cannot be established on (ZF-AC);
because of Proposition \ref{p5}.
\erem

Finally, an interesting question to be posed 
is that of such inclusion chains
being retainable beyond the countable case.
Unfortunately, this is not in general possible;
see
Isac \cite{isac-1983}
for details.

\section{Gauge variational principles}
\setcounter{equation}{0}

Let $(X,\le)$ be a quasi-ordered structure.
By a {\it pseudometric} over $X$ we shall
mean any map $d:X\times X \to R_+$. If, in addition, $d$ is
{\it triangular} [$d(x,z)\le d(x,y)+d(y,z), \forall x,y,z\in X$],
{\it symmetric} [$d(x,y)=d(y,x), \forall x,y\in X$] and
{\it reflexive} [$d(x,x)=0, \forall x\in X$],
we say that it is a {\it semimetric} over $X$.
Suppose that we fixed such an object;
the triple $(X;\le;d)$ will be then referred to as a 
{\it quasi-ordered semimetric space}.
The sequential convergence $(\dtends)$ attached to $d$ means:
the sequence $(x_n)$ in $X$, $d$-{\it converges} to $x\in X$ 
(and we write: $x_n \dtends x$), 
iff $d(x_n,x)\to 0$ as $n\to \oo$;
and reads: $x$ is a $d$-limit of $(x_n)$.
When $x$ is generic in this convention we say that $(x_n)$ is
$d$-{\it convergent}.
Further, the $d$-Cauchy property 
of a sequence $(x_n)$ in $X$ means:
$d(x_m,x_n)\to 0$ as $m,n\to \oo$.
By the imposed upon $d$ properties, 
each $d$-convergent sequence is 
$d$-Cauchy too; the reciprocal is not in general valid.
\sk

{\bf (A)}
Let $D=(d_i;i\ge 0)$  be a denumerable 
family of semimetrics on $X$; supposed to be
{\it sufficient} [$d_i(x,y)=0$, $\forall i\ge 0$, implies $x=y$];
in this case, the couple $(X;\le;D)$
will be termed a {\it quasi-ordered gauge space}.
We say that the sequence $(x_n; n\ge 0)$ in $X$, 
$D$-{\it converges} to $x\in X$ 
(and we write $x_n\Dtends x$), 
when it $d_i$-converges to $x$, for each $i\ge 0$. 
Likewise, the sequence $(x_n)$ in $X$ is called $D$-{\it Cauchy}, 
when it is $d_i$-Cauchy, for each $i\ge 0$.
By the remark above, any (ascending) $D$-convergent sequence is 
(ascending) $D$-Cauchy. 
If the reciprocal holds -- for ascending sequences -- 
then $(X;\le;D)$ is termed {\it complete}. 
Call the subset $Z$ of $X$, $(\le,D)$-{\it closed} when the $D$-limit
of each ascending sequence in $Z$ belongs to $Z$.
In particular, we say that $(\le)$ is {\it $D$-self-closed} provided
$X(x,\le)$ is $(\le,D)$-closed, for each $x\in X$; 
or, equivalently: the $D$-limit of each ascending sequence
is an upper bound of it (modulo $(\le)$).
By definition, the property
[$(X;\le;D)$ is complete and $(\le)$ is $D$-self-closed]
will be referred to as: $(X;\le;D)$ is {\it strongly complete}.

Having these precise, let us introduce a lot of
(topological type) subclasses (=subsets) of $\calf(X,\Rpoo)$ as
\bit
\item[(L1)]
$\vphi$ is {\it descending $(\le,D)$-lsc}:\
$(x_n)$=ascending, $x_n\Dtends x$, \\
($x_n\le x$, $\forall n$)\
and $(\vphi(x_n))$=descending,\
imply  [$\vphi(x_n)\ge \vphi(x)$, $\forall n$]
\item[(L2)] 
$\vphi$ is {\it descending $D$-lsc}:\\
$x_n\Dtends x$, and $(\vphi(x_n))$=descending, imply 
[$\vphi(x_n)\ge \vphi(x)$, $\forall n$]
\item[(L3)] 
$\vphi$ is {\it $D$-lsc}:\
$\liminf_n \vphi(x_n)\ge \vphi(x)$,\ whenever $x_n\Dtends x$.
\eit
Note that (L3) $\limpl$ (L2) $\limpl$ (L1); 
we do not give details.
Further, again via (c05), these 
give corresponding subclasses (=subsets) ([Lk]; $k\in \{1,2,3\}$),
in the gauge functions class 
over $\calf(X,\Rpoo)$; hence, by the above,
[L3] $\limpl$ [L2] $\limpl$ [L1].

The following (quasi-order) "multiple" 
gauge variational principle is now entering into 
our discussion:

\btheorem \label{t4}
Assume that $(X;\le;D)$ is strongly complete;
and let $\Phi=(\vphi_i; i\ge 0)$ be a gauge function over
$\calf(X,\Rpoo)$, fulfilling
\bit
\item[(d01)]
$\Phi$ is proper:
$\Dom(\Phi):=\cap\{\Dom(\vphi_i); i\ge 0]\}\ne \es$
\item[(d02)]
$\Phi$ belongs to the subclass [Pj], for some $j\in \{3,4\}$
\item[(d03)]
$\Phi$ belongs to the subclass [Lk], for some $k\in \{1,2,3\}$.
\eit
Then, for each $u\in \Dom(\Phi)$ there exists 
$v\in \Dom(\Phi)$ with 
\beq \label{401}
u\le v,  d_i(u,v)\le \vphi_i(u)-\vphi_i(v),\ \ \forall i\ge 0 
\eeq
\beq \label{402}
\forall x\in X(v,\le)\sm \{v\}, \exists i=i(x):\
d_i(v,x)> \vphi_i(v)-\vphi_i(x).
\eeq
\etheorem

By definition, this result will be 
written as (EVPg-Pj-Lk).
Note that, by the inclusions above,
\beq \label{403}
\mbox{
(EVPg-Pj-L1) $\limpl$ (EVPg-Pj-L2)
$\limpl$ (EVPg-Pj-L3),\ $\forall j\in \{3,4\}$.
}
\eeq
On the other hand, 
by the argument in Lemma \ref{le1},
\beq \label{404}
\mbox{
(EVPg-P3-Lk) $\lequi$ (EVPg-P4-Lk),\
$\forall k\in \{1,2,3\}$.
}
\eeq

As in Section 3, the obtained relations cannot assure us that 
these principles are deductible in (ZF-AC+DC).
This, however, holds; as results from

\bprop \label{p8}
We have (in (ZF-AC)) (DC) $\limpl$ (EVPg-P3-L1);
hence (by the above) 
(DC) $\limpl$ (EVPg-Pj-Lk), 
for all $j\in \{3,4\}$
and all $k\in \{1,2,3\}$.
\eprop

\bproof
Let $(\preceq)$ stand for the quasi-order (over $X$):
\bit
\item[]
$x\preceq y$ iff $x\le y$\ and\ 
$(d_i(x,y)+\vphi_i(y)\le \vphi_i(x)$, $\forall i$).
\eit
Clearly, $(\preceq)$ is antisymmetric
[hence, an order]
on $\Dom(\Phi)$; 
so, it remains as such over its subset $X[u]:=X(u,\preceq)$.
We claim that conditions of (BBg) 
(i.e.: the gauge ordering principle in
Turinici \cite{turinici-1982})
are fulfilled over
$(X[u];\preceq;\Phi)$. 
Clearly, $\Phi$ satisfies (c04); and 
(thanks to [P3]), it belongs to the subclass [P5] 
relative to $(X[u],\preceq)$;
so, it remains to show that 
(a01) holds over $(X[u],\preceq)$.
Let $(x_n)$  be an $(\preceq)$-ascending sequence in $X[u]$:
\bit
\item[(d04)]
$x_n\le x_m$\ and  
[$d_i(x_n,x_m)\le \vphi_i(x_n)-\vphi_i(x_m)$,\ $\forall i\ge 0$],\ if $n\le m$.
\eit
By [P5] (see above), it follows that, for each $i\ge 0$, 
the sequence $(\vphi_i(x_n))$ is descending and 
bounded from below; hence a Cauchy one.
This, along with (d04), 
tells us that $(x_n)$ 
is a $(\le)$-ascending $D$-Cauchy sequence in $X[u]$; 
wherefrom (by  completeness), 
there must be some $y\in X$ with $x_n\Dtends y$.
Note that, by the self-closedness property, $x_n\le y$, $\forall n$;
and this, via [L1], yields 
[$\vphi_i(x_n)\ge \vphi_i(y)$, $\forall i$, $\forall n$].
For each pair $(i,n)$, we have
$$ \barr{l}
d_i(x_n,y)\le d_i(x_n,x_m)+d_i(x_m,y)\le 
\vphi_i(x_n)-\vphi_i(x_m)+d_i(x_m,y)\le \\
\vphi_i(x_n)-\vphi_i(y)+d_i(x_m,y),\ \ \forall m\ge n.
\earr
$$
Passing to limit as $m\to \oo$ one derives 
(in combination a previous fact)
$$ 
\mbox{
($\forall n$):\ 
$x_n\le y, [d_i(x_n,y)\le \vphi_i(x_n)-\vphi_i(y),\  \forall i]$;\
hence,  $x_n\preceq y$. 
}
$$
This firstly shows that $y\in X[u]$;
and secondly, that $y$ is an upper bound (modulo $(\preceq)$) in $X[u]$ of $(x_n)$. 
Summing up, $(X[u],\preceq)$ is sequentially inductive;
as claimed.
From (BBg) it follows that, for the starting 
$u\in X[u]$, there exists $v\in X[u]$ with

{\bf h)} $u\preceq v$;\   
{\bf hh)} $v\preceq x\in X[u]$ $\limpl$ 
[$\vphi_i(v)=\vphi(x)$, $\forall i$].

\n
The former of these is just (\ref{401}). And the latter one gives at once
(\ref{402}). In fact, let $y\in X$ be such that 
$v\le y$, [$d_i(v,y)\le \vphi_i(v)-\vphi_i(y)$, $\forall i$]. 
As a consequence, 
$v\preceq y\in X[u]$; 
so that (by {\bf hh)} above) 
$\vphi_i(v)=\vphi_i(y)$, $\forall i$.
This, by the working hypothesis above yields 
[$d_i(v,y)=0$, $\forall i$];
so that (as $D$ is sufficient) $v=y$. 
The proof is complete.
\eproof

A basic particular case of our developments corresponds to
$(\le)=X\times X$ (=the {\it trivial} quasi-order on $X$).
Then, the subclasses (L1) and (L2) are identical;
hence, so are the gauge subclasses [L1] and [L2].
By Theorem \ref{t4}  we get the (amorphous)
"multiple" gauge variational principle:

\btheorem \label{t5}
Assume that $(X,D)$ is complete;
and let $\Phi=(\vphi_i; i\ge 0)$ be a gauge function over 
$\calf(X,\Rpoo)$, fulfilling
(d01), (d02), and 
\bit
\item[(d05)]
$\Phi$ is in the subclass [Lk], for some $k\in \{2,3\}$.
\eit
Then, for each  $u\in \Dom(\Phi)$ there exists
$v=v(u)\in \Dom(\Phi)$ with
\beq \label{405}
d_i(u,v)\le \vphi_i(u)-\vphi_i(v),\ \  \forall i\ge 0
\eeq
\beq \label{406}
\forall x\in X\sm \{v\}, \exists i=i(x):\
d_i(v,x)> \vphi_i(v)-\vphi_i(x).
\eeq
\etheorem

Finally, a particular case of these facts is that of the gauge function 
$\Phi=(\vphi_i; i\ge 0)$
having all components with finite values;
hence, $\Dom(\Phi)=X$. 
Then, from Theorem \ref{t5}, we get
the (finitary) "multiple" gauge variational principle:

\btheorem \label{t6}
Assume that $(X,D)$ is complete; and let 
$\Phi=(\vphi_i; i\ge 0)$ be a gauge function over
$\calf(X,R)$, fulfilling (d02) and (d05).
Then, for each $u\in X$, there exists $v=v(u)\in X$
with the properties (\ref{405}) and (\ref{406}).
\etheorem

\brem \label{r3}
\rm

As shown in (\ref{404}),
the results (EVPg-P3-L3) and (EVPg-P4-L3)
are equivalent in (ZF-AC+DC) and both
of these are deductible from (DC) in (ZF-AC).
Note that, in the context of Theorem \ref{t6},
the former of them is just the
gauge variational principle in 
Turinici \cite{turinici-1982}:
while the latter is the 
gauge variational statement in 
Bae et al \cite{bae-cho-kim-2011}.
\erem

{\bf (C)}
Now, assume that, in these results, 
$D$ and $\Phi$ are constant sequences.
We have three cases to consider.

{\bf I)}
Let $(X;\le;d)$ be a quasi-ordered metric space.
Given a function $\vphi:X\to \Rpoo$, 
call it {\it descending $(\le,d)$-lsc} provided
(L1) holds, with $D=\{d\}$.
As a direct consequence of Theorem \ref{t4}, we have
the monotone variational principle in
Turinici \cite{turinici-1990}
(in short: (EVP-m)).

\bcor \label{c1}
Assume that $(X;\le;d)$ is strongly complete;
and that the function $\vphi:X\to \Rpoo$ 
is proper, bounded from below and descending 
$(\le,d)$-lsc.
Then,  for each $u\in \Dom(\vphi)$ there exists 
$v\in \Dom(\vphi)$ with 
\beq \label{407}
u\le v,  d(u,v)\le \vphi(u)-\vphi(v) 
\eeq
\beq \label{408}
\forall x\in X(v,\le)\sm \{v\}:\
d(v,x)> \vphi_i(v)-\vphi_i(x).
\eeq
\ecor

The motivation of our terminology
comes from the fact that 
(L1) holds under (a02)
(and the $d$-self-closedness of $(\le)$).
This determines us to consider 
the components of Theorem \ref{t4}, 
as gauge versions of (EVP-m).

{\bf II)}
Let $(X,d)$ be a metric space.
Call the function $\vphi:X\to \Rpoo$,
{\it descending $d$-lsc} provided (L2)
holds, with $D=\{d\}$.
As a direct consequence of Theorem \ref{t5}, 
we get the (descending) 
Ekeland's variational principle \cite{ekeland-1979}
(in short: (EVP)).

\bcor \label{c2}
Assume that $(X,d)$ is complete;
and that the function $\vphi:X\to \Rpoo$ 
is proper, bounded from below and descending $d$-lsc.
Then,  for each $u\in \Dom(\vphi)$ there exists 
$v\in \Dom(\vphi)$ with 
\beq \label{409}
d(u,v)\le \vphi(u)-\vphi(v) 
\eeq
\beq \label{410}
\forall x\in X\sm \{v\}:\
d(v,x)> \vphi(v)-\vphi(x).
\eeq
\ecor

By this relationship, it is natural that 
the components of Theorem \ref{t5} be taken
as gauge  versions of (EVP);
we shall denote them as (EVPg).
On the other hand, (EVP-m) includes (EVP):
just take $(\le)$ as the trivial quasi-order on $X$.

{\bf III)}
Let $(X,d)$ be a metric space.
Given a function $\vphi:X\to R$, 
call it {\it $d$-lsc} provided (L3) holds,
with $D=\{d\}$.
As a direct consequence of Theorem \ref{t7},
we get the finitary 
Ekeland's variational principle
(in short: (EVP-f)).

\bcor \label{c3}
Assume that $(X,d)$ is complete;
and that the function $\vphi:X\to R$ 
is bounded from below and $d$-lsc.
Then,  for each $u\in X$ there exists 
$v\in X$ with 
the properties (\ref{409}) and (\ref{410}).
\ecor

As before, this relationship determines us to 
consider the components  of Theorem \ref{t6}
as gauges version of (EVP-f).
Moreover, (EVP) includes (EVP-f) in a trivial way. 
\sk

Finally, we have to stress that many other 
maximal/variational statements have gauge versions.
Further aspects will be delineated elsewhere.

\section{(EVPdLc) implies (DC)}
\setcounter{equation}{0}

As a consequence of the previous developments,
all results in Section 3 and Section 4 are in the 
logical segment between (DC) and (EVP-f).
So, it is natural to determine the "logical" ecart
between these extreme terms.
The natural setting for solving this problem 
is (ZF-AC)(=the {\it reduced} Zermelo-Fraenkel system).
\sk

Let $X$ be a nonempty set; 
and $(\le)$ be an order on it. We say that $(\le)$
has  the {\it inf-lattice} property, provided: 
$x\wedge y:=\inf(x,y)$ exists, for all $x,y\in X$.
Further, we say that $z\in X$ is a $(\le)$-{\it maximal} element if 
$X(z,\le)=\{z\}$; the class of all these points will be
denoted as $\max(X,\le)$. 
In this case, $(\le)$ is called a {\it Zorn order} when 
$\max(X,\le)$ is nonempty and {\it cofinal} in $X$
[for each $u\in X$ there exists a $(\le)$-maximal 
$v\in X$ with $u\le v$].
Further aspects are to be described in a metric setting. 
Let $d:X\times X\to R_+$ be a metric over $X$;
and $\vphi:X\to R_+$ be some function.
Then, the natural choice for $(\le)$ above is
\bit
\item[]
$x\le_{(d,\vphi)} y$ iff $d(x,y)\le \vphi(x)-\vphi(y)$;
\eit
referred to as the
Br{\o}ndsted order \cite{brondsted-1976}
attached to $(d,\vphi)$. 
Denote
$X(x,\rho)=\{u\in X; d(x,u)< \rho\}$, $x\in X$, $\rho> 0$
[the open sphere with center $x$ and radius $\rho$].
Call the ambient metric space  $(X,d)$, {\it discrete} when for
each $x\in X$ there exists $\rho=\rho(x)> 0$ such that 
$X(x,\rho)=\{x\}$. Note that, under such an assumption, 
any function $\psi:X\to R$ is continuous over $X$.
However, the Lipschitz property 
($|\psi(x)-\psi(y)|\le L d(x,y)$, $x,y\in X$, for some $L> 0$)
cannot be assured, in general.
\sk

Now, the variational principle below enters into our discussion:

\btheorem \label{t7}
Let  the metric space $(X,d)$ and the function $\vphi:X\to R_+$
satisfy
\bit
\item[(e01)]
$(X,d)$ is discrete bounded and complete
\item[(e02)]
$(\le_{(d,\vphi)})$ has the inf-lattice property
\item[(e03)]
$\vphi$ is $d$-nonexpansive and $\vphi(X)$ is countable.
\eit
Then, $(\le_{(d,\vphi)})$ is a Zorn order.
\etheorem

We shall refer to it as: the discrete Lipschitz countable version of EVP 
(in short: (EVPdLc)).
Clearly, (EVP-f) $\limpl$ (EVPdLc). 
The remarkable fact to be added is that
this last principle yields (DC); so, it 
completes the circle between all these.

\bprop \label{p9}
We have (in (ZF-AC)) (EVPdLc) $\limpl$ (DC).
So, the gauge (standard)  ordering/variational principles
in Section 3 -- Section 5
are all equivalent with (DC);
hence, mutually equivalent.
\eprop

For a complete proof, see 
Turinici \cite{turinici-2011}.
In particular, when the specific assumptions 
(e02) and (e03) (the second half)
are ignored in Theorem \ref{t7}, 
Proposition \ref{p9}
is comparable with a related statement in 
Brunner \cite{brunner-1987}. 
Further aspects may be found in 
Dodu and Morillon \cite{dodu-morillon-1999};
see also 
Schechter \cite[Ch 19, Sect 19.53]{schechter-1997}.


\end{document}